\documentclass[12pt]{article}

\usepackage{multirow}

\usepackage{enumitem}
\usepackage{subfig}

\usepackage{lipsum} 
\usepackage{wrapfig}
\usepackage{graphicx}
\usepackage{ifthen}
\newboolean{PrintVersion}
\setboolean{PrintVersion}{false} 
\usepackage{amsmath,amssymb}
\usepackage{amsthm}
\usepackage{mathtools}
\usepackage{graphicx}
\usepackage{caption}
\usepackage{bbm}
\usepackage{color}
\usepackage{dcolumn}
\usepackage{bm}
\usepackage{float}
\usepackage[pagebackref=true]{hyperref}
\renewcommand*\backref[1]{\ifx#1\relax \else (Cited on #1) \fi}\usepackage{apacite}
\definecolor{background-color}{gray}{0.98}
\definecolor{steelblue}{rgb}{0.27, 0.51, 0.71}
\definecolor{brickred}{rgb}{0.8, 0.25, 0.33}
\definecolor{bluegray}{rgb}{0.4, 0.6, 0.8}
\definecolor{amethyst}{rgb}{0.6, 0.4, 0.8}

\hypersetup{
    plainpages=false,       
    unicode=false,          
    pdftoolbar=true,        
    pdfmenubar=true,        
    pdffitwindow=false,     
    pdfstartview={FitH},    
    pdfnewwindow=true,      
    colorlinks=true,        
    linkcolor=steelblue,         
    citecolor=bluegray,        
    filecolor=magenta,      
    urlcolor=cyan           
}
\ifthenelse{\boolean{PrintVersion}}{   
\hypersetup{	
    citecolor=black,%
    filecolor=black,%
    linkcolor=black,%
    urlcolor=black}
}{} 

\theoremstyle{plain}
\newtheorem{theorem}{Theorem}[section]

\newtheorem{corollary}[theorem]{Corollary}

\newtheorem{proposition}[theorem]{Proposition}

\theoremstyle{definition}

\newenvironment{classification}
  {\noindent \textit{AMS classification:}}
  
\newenvironment{keywords}
  {\noindent \textit{Keywords:}}
  
\newcommand{\ve}[1]{\mathbf{#1}}           
\newcommand{\sv}[1]{\boldsymbol{#1}}   
\newcommand{\m}[1]{\mathbf{#1}}               
\newcommand{\tr}[1]{{#1}^{\mkern-1.5mu\mathsf{T}}}              
\newcommand{\inv}[1]{{#1}^{\mkern-1.5mu{-1}}}              

\newcommand*{\intersect}{\cap}
\newcommand*{\union}{\cup}

\newcommand*{\comp}[1]{\overline{#1}} 
\newcommand{\cardinality}[1]{|{#1}|}

\newcommand*{\powerset}[1]{\mathcal{P}({#1})} 

\newcommand{\family}[1]{{\cal #1}}

\let\emptyset\varnothing

\newcommand{\Nset}[1]{\mathcal{N}_{#1}}

\newcommand{\mCliqueSet}[1]{\{#1_1, \ldots, #1_m \}}  
\newcommand{\support}[1]{Supp(#1)}

\newcommand{\intersectingFamily}[1]{\mathcal{I}}
\newcommand{\pathIntersectingFamily}[1]{\mathcal{I}_p}
\newcommand{\field}[1]{\mathbb{#1}}

\newcommand{\NaturalsZero}{\field{N}_0}

\newcommand{\imply}{\Longrightarrow}

\title{How many cliques can a clique cover cover?}

\author{
Pavel Shuldiner \& R. Wayne Oldford\\
  University of Waterloo
}

\date{}
\begin{document}
\maketitle
\begin{abstract}
This work examines the problem of clique enumeration on a graph by exploiting its clique covers. The principle of inclusion/exclusion is applied to determine the number of cliques of size $r$ in the graph union of a set $\family{C} = \{c_1, \ldots, c_m\}$ of $m$ cliques. This leads to a deeper examination of the sets involved and to an orbit partition, $\Gamma$, of the power set $\powerset{\Nset{m}}$ of $\Nset{m} = \{1, \ldots, m\}$. Applied to the cliques, this partition gives insight into clique enumeration and yields new results on cliques within a clique cover, including expressions for the number of cliques of size $r$ as well as generating functions for the cliques on these graphs. The quotient graph modulo this partition provides a succinct representation to determine cliques and maximal cliques in the graph union. The partition also provides a natural and powerful framework for related problems, such as the enumeration of induced connected components, by drawing upon a connection to extremal set theory through intersecting sets. 
\end{abstract}
\section{Introduction}
\label{section:intro}
%
For any graph $G = (V, E)$ and node subset $H \subseteq V$, the induced subgraph $G[H]$ has nodes $H$ and those edges in $E$ whose endpoints lie in $H$.  A clique of size $r$ is induced whenever $G[H]$ is a complete graph on $r$ nodes.  Allowing trivial cliques (i.e., $r =1$ or $r=2$), a collection of cliques  $\family{C} =\mCliqueSet{c}$ can always be found (for some $m$) which \textit{covers} the graph $G$ -- in the sense that the graph union, $G[c_1] \union G[c_2] \union \cdots \union G[c_m]$, of the induced subgraphs has the same vertex set as $G$.  

Such a collection is called a \textit{vertex clique cover} of $G$.
When its cliques are also non-intersecting  (i.e., $c_i \intersect c_j = \emptyset ~ \forall i \ne j$), then the collection will be called a \textit{clique cover partition}, so as to clearly identify this special case.

A collection of cliques whose graph union contains all edges in $G$ is called an \textit{edge clique cover} \cite<e.g., see>{edgeCliques85}.
In what follows, interest lies in counting the number of cliques, of any specified size $r$, formed by the graph union over   \textit{any} of these clique covers, indeed over any \textit{collection of induced cliques} of $G$.  

Suppose the graph $G$ has $n$ nodes numbered 1 to $n$, so that the power set, $\powerset{\Nset{n}}$, of $\Nset{n} = \left\{1, \ldots, n \right\}$ identifies, by node indices, all possible induced subgraphs of $G$.  
For index set $i \subseteq \Nset{n}$, provided $G$ is understood, the induced subgraph $G[i]$ may be  more simply denoted by its index set $i$.
A collection of cliques, then, is denoted by a \textit{family of sets} $\family{C} = \left\{c_1, \ldots, c_m \right\} \subseteq \powerset{\Nset{n}}$, provided each $c_j \in \family{C}$ identifies a clique induced in $G$.

\begin{figure}[h]
\begin{center}
\includegraphics[width = 0.45\textwidth]{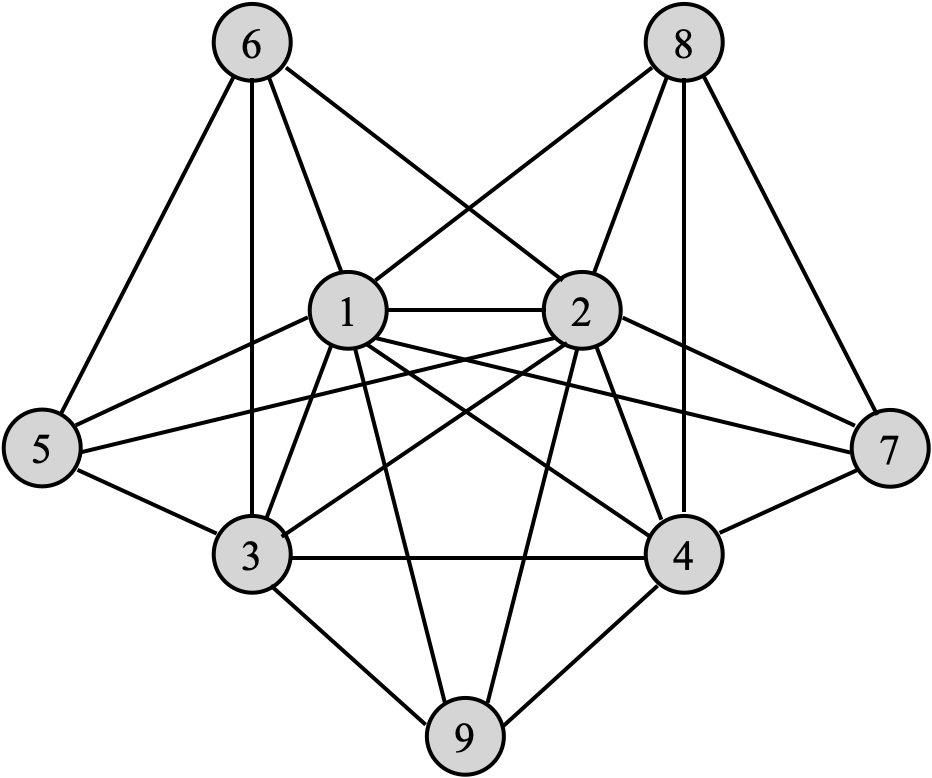}
\end{center}
  \caption{The graph union of $\family{C} = \left\{A, B, C\right\}$ with $A := \{1, 2, 3, 5, 6\}, B:= \{1,2,4,7,8\}$ and $C:=\{1,2,3,4,9\}$.  How many cliques are there of size $r=1, 2, 3, \ldots$?}
\label{fig:clique_collection}
\end{figure}

For example, suppose $n \ge 9$ and $G$ contains three cliques $A := \{1, 2, 3, 5, 6\}, B:= \{1,2,4,7,8\}$ and $C:=\{1,2,3,4,9\}$, each of size $5$.  Then $\family{C} = \left\{A, B, C\right\}$ is a collection of three size $5$ cliques, being a \textit{vertex clique cover} only if $n = 9$ (and not if $n > 9$).  Its graph union is shown in Figure \ref{fig:clique_collection}.
It may, or may not, also be an \textit{edge clique cover}, depending on whether, or not, the union contains \textit{all} edges of $G$.  It is \textit{not} a  \textit{clique cover partition} because the intersection of at least one pair of $A$, $B$, and $C$ is non-null (here all pairs intersect).

Our interest lies in determining the number of cliques of any size $r$ in the union.  When $r=5$, there are exactly three $5$-cliques, namely $A$, $B$ and $C$.
Consulting Figure \ref{fig:clique_collection}, there are no cliques in the union of size $r \ge 6$, though this need not be true in general -- e.g., were the $3$-clique $D:=\{4,5,6\}$ added to the collection $\family{C}$, the $6$-clique $\{1, 2, 3, 4, 5, 6\}$ would arise. 
For $r <5$, the intersections of the cliques in $\family{C}$ must also be considered.
If all intersections are null, then $\family{C}$ would be \textit{clique cover partition} of its union, and the number of cliques of size $r \le 5$ would simply be the sum of the number of $r$-cliques within each clique of $\family{C}$. But that is not the case here -- e.g., the $3$-clique $\left\{1, 2, 3 \right\}$ appears in both $A$ and $C$ -- so care is needed to avoid overcounting. 
Careful examination of Figure \ref{fig:clique_collection} will yield 15 cliques of size 4, 28 of size 3, and 24 of size 2.

Given a collection $\family{C} = \{c_1, \ldots, c_m\}$ of $m$ cliques on the same graph, an expression for the number of cliques of size $r$ in the graph union $G = \union_{i=1}^m c_i$ can be had by applying the the principle of inclusion/exclusion.  This is done in Section \ref{section:PIE_for_cliques}. 
 
A richer approach is to first form a partition of $G = \union_{i=1}^m c_i$ based on index sets $J$, now consisting of the indices from $\{1, \ldots, m\}$ which identify the cliques in the collection (or cover) $\family{C}$.  That is, each partition cell is identified with one set $J \in \powerset{\Nset{m}}$; the set of graph node indices in cell $J$ will be denoted $\Gamma_J \in \powerset{\Nset{n}}$ and the partition called a $\Gamma$-partition.  
The cardinality of $\Gamma_J$ will be denoted $\gamma_J$.
This is the primary approach introduced and explored in this paper.

Subgraphs $H$ of $G = \union_{i=1}^m c_i$ will have nodes appearing in some cells $\Gamma_J$ (for some $J \in \powerset{\Nset{m}}$) and not in others.  The clique index cells $J$ whose $\Gamma_J$ contain nodes in the subgraph $H$ will be called the \textit{support} of $H$ and the tuple containing the count of nodes of $H$ in each $\Gamma_J$ its \textit{signature}.  Whether $H$ is connected, or is a clique of size $r$, or forms a maximal clique, can be determined using characteristics of its support and/or signature.  This is shown by connecting these concepts to intersecting sets and \textit{intersecting families} of sets \cite<e.g., see>{meyerowitz1995maximal}.

The $\Gamma$-partition is itself an \textit{orbit} partition and hence an \textit{equitable} partition.  The quotient graph, $G/\Gamma$, which results compresses and contains all information needed to determining connected subgraphs, cliques, and maximal cliques on $G$.

Section \ref{section:partition-example} introduces and illustrates these concepts using the three clique collection example of Figure \ref{fig:clique_collection}.  Section \ref{section:partitional-framework} then provides a more general treatment with formal definitions and proved results. The general $\Gamma$-partition is derived in Section \ref{section:partition-general} for any collection of subsets of $\Nset{m}$ and applied to clique collections in Section \ref{section:cliquepartition} where.  It is shown to be an orbit partition in Section \ref{section:generalquotient} and its quotient graph defined.  Support and signatures are formally defined in Section \ref{section:types} and used to define different types of isomorphic graphs.  Section \ref{section:signatures} establishes some counting results on signatures as does Section \ref{section:connected-subgraphs} as they relate to subgraph connectedness.  Section \ref{section:connected-subgraphs} ends with a generating function for the number of induced connected subgraphs of size $k$.  Section \ref{section:supportandcliques} shows $H$ induces a clique if, and only if, its support is an intersecting family; Theorem \ref{thm:maximal_clique} provides the conditions for the clique to be maximal.  Section \ref{quotientandmaximalcliques} shows how the quotient graph, $G/\Gamma$, can be used to directly determine cliques and maximal cliques in the original graph union $G$ and ends with some minor results on the number of maximal cliques and the clique number of $G$.

Section \ref{section:cliquecounting} uses the framework of Section \ref{section:partitional-framework} to finally get down to counting cliques.  
Results include expressions for the number of cliques containing any particular clique $H$,  the number of cliques of size $r$, and, in Theorem \ref{thm:clique_gen_fn}, a generating function for clique counts in the graph union of $m$ cliques.  Theorem  \ref{thm:clique_gen_fn} is then applied to give a new expression for the number of $r$-cliques and for the number of edges induced by a collection of $m$ cliques of size $r$.
The section ends with an application of the ``hand-shaking lemma'' to yield an expression for the number of edges induced by any collection of cliques.
%
%

The paper ends with a brief summary discussion as Section \ref{section:discussion}.

\section{Counting by inclusion/exclusion}
\label{section:PIE_for_cliques}
As the example of Section \ref{section:intro} suggests, the key to clique counting over the graph union of a collection of cliques will be identifying the intersection of the various index sets.  
Unsurprisingly, then, our first approach to enumerating cliques makes use of the \textit{Principle of Inclusion/Exclusion}.

This yields the following result for the count of the number of $r$-cliques in the union of an arbitrary collection of cliques.
\begin{proposition}
\label{prop:PIE_for_cliques}
	Let $\mathcal{C} = \{c_1, \ldots, c_m\}$ be a collection of cliques. The total number of $r-$cliques that are induced by $\{c_1, \ldots, c_m\}$ is
	$$\sum_{J:\emptyset \neq J \subseteq\{1,\ldots, m\}}
(-1)^{|J|+1}\binom{I_J}{r},$$
where $I_J := |\bigcap_{j\in J}c_j|$.
\end{proposition}

\begin{proof}
We count the number of $r-$cliques that are induced by at least one of the cliques in $\mathcal{C}$. Let $\binom{c_j}{r} := \{\{v_1,\ldots, v_r\} \subseteq c_j: v_1\neq  \dotsb \neq v_r\}$ denote the set of $r-$cliques induced by the clique $c_j$.
We will prove that for any nonempty $J\subseteq \{1,\ldots, m\}$,
$$\left|\bigcap_{j\in J} \dbinom{c_j}{r}\right| 
= \binom{I_J}{r},$$
by showing that 
$$\bigcap_{j\in J} \dbinom{c_j}{r} = \dbinom{\bigcap_{j\in J}c_j}{r}.$$ 
If $\{v_1,\ldots, v_r\} \in \bigcap_{j\in J} \binom{c_j}{r}$, then $\{v_1,\ldots, v_r\}  \subset c_j$ for all $j\in J$ and so
$$\{v_1,\ldots, v_r\}  \in \binom{\bigcap_{j\in J}c_j}{r}.$$
Conversely, if $\{v_1,\ldots, v_r\} \in \dbinom{\bigcap_{j\in J}c_j}{r}$ then $\{v_1,\ldots, v_r\} \subset c_j$ for all $j\in J$. Therefore,
$$\{v_1,\ldots, v_r\} \in c_j,$$
for all $j\in J$ and the claim follows. 

Therefore, the total number of $r-$cliques within $A$ is $\left|\bigcup_{j\in \{1,\ldots, m\}} \dbinom{c_j}{r}\right|$. By Principle of Inclusion Exclusion \cite<e.g., see>[p. 112]{wilf2005generatingfunctionology}, 

$$\left|\bigcup_{j\in J} \dbinom{c_j}{r}\right| = \sum_{\emptyset \neq J \subseteq\{1,\ldots, m\}}
(-1)^{|J|+1}\left|\bigcap_{j\in J} \dbinom{c_j}{r}\right|
= \sum_{\emptyset \neq J \subseteq\{1,\ldots, m\}}
(-1)^{|J|+1}\binom{I_J}{r},$$
as needed to be shown.
\end{proof}
This leads to an expression for the \textit{total} number of typically interesting cliques (i.e., $r \ge 3$; non-trivial: no single edge, no single vertex, cliques):
\begin{corollary}
	\label{cor:total_num_induced_cliques}
	Let $\mCliqueSet{c}$ be a collection of cliques. The total number of non-trivial cliques that are contained within $\mCliqueSet{c}$ is 
	$$\sum_{J:\emptyset \neq J \subseteq\{1,\ldots, m\}}
(-1)^{|J|+1}\left(2^{I_J}-\binom{I_J}{2} \right) - \left|\bigcup_{j=1}^{m}c_j \right| -1 .$$
\end{corollary}
\begin{proof}
	By Proposition $\ref{prop:PIE_for_cliques}$, the total number of cliques is given by
	\begin{align*}
		\sum_{r = 0}^{\infty}\sum_{J:\emptyset \neq J \subseteq\{1,\ldots, m\}}
			(-1)^{|J|+1}\binom{I_J}{r} &= 
				\sum_{J:\emptyset \neq J \subseteq\{1,\ldots, m\}}
						(-1)^{|J|+1}\sum_{r = 0}^{\infty}\binom{I_J}{r} \\
						&= 
							\sum_{J:\emptyset \neq J \subseteq\{1,\ldots, m\}}
								(-1)^{|J|+1}\sum_{r = 0}^{\infty}2^{I_J},
	\end{align*}
	by the Binomial Theorem. Now, since there is only one $0-$clique on a set of nodes, the 1-cliques correspond to the $\left|\bigcup_{j=1}^{m}c_j\right|$ vertices and 2-cliques is the number of edges, 
	$$\sum_{J:\emptyset \neq J \subseteq\{1,\ldots, m\}}
								(-1)^{|J|+1}\sum_{r = 0}^{\infty}2^{I_J} = \sum_{J:\emptyset \neq J \subseteq\{1,\ldots, m\}}
(-1)^{|J|+1}\left(2^{I_J}-\binom{I_J}{2} \right) - \left|\bigcup_{j=1}^{m}c_j \right| -1.$$
\end{proof}

For example, let $r = 3$ and let $\mathcal{C}= \{c_1, c_2\}$ be a collection of the two triangles $c_1 = \{1, 2, 3\}$ and $c_2 = \{2, 3, 4\}$.  If $e_j$ is the number of edges induced by triangle $j$, then the total number of edges in the collection is given by
$$e_1 + e_2 - \binom{|c_1\bigcap c_2|}{2} = \binom{|c_1|}{2} + \binom{|c_2|}{2} - \binom{2}{2} = 3+3-1 = 5$$
since $\dbinom{|c_1\bigcap c_2|}{2}$ is the number of edges common to both $c_1$ and $c_2$ (one edge for every 2 vertices). 

\section{A partition framework}
\label{section:partition-example}
\begin{figure}[htbp]
\begin{center}
\begin{tabular}{cr}
\begin{tabular}{c}
\includegraphics[width = 0.35\textwidth]{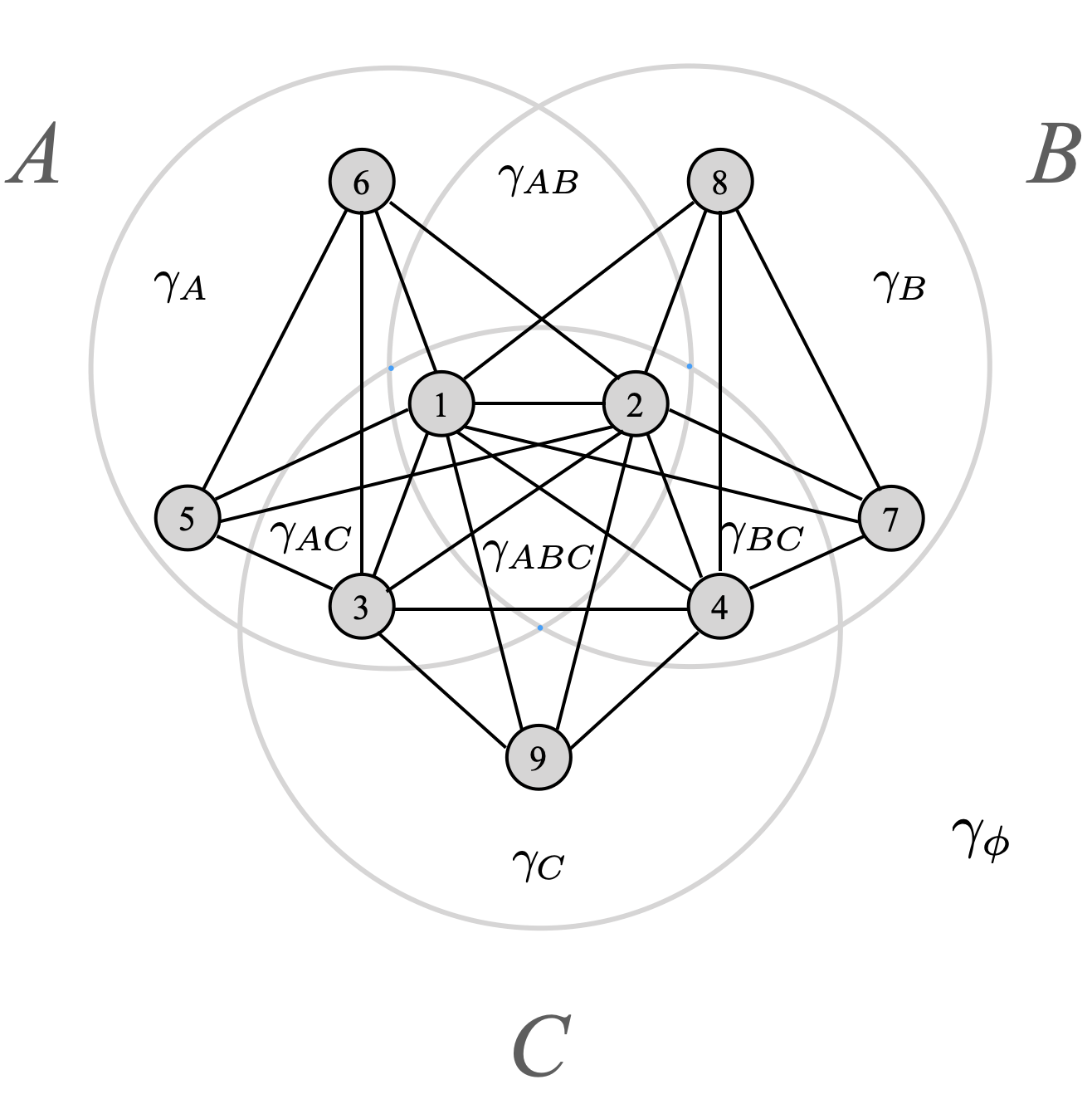}  
\end{tabular}
&
\begin{tabular}{rcl}
 $ \Gamma_{A} = \{5,6\} $ & $\implies$ &  $ \gamma_{A} = 2$ \\
 $ \Gamma_{B} = \{7,8\} $ & $\implies$ &  $ \gamma_{B} = 2$ \\
 $ \Gamma_{C} = \{9\}  $ & $\implies$ &  $ \gamma_{C} = 1$ \\
 $ \Gamma_{AB} =  \emptyset  $ & $\implies$ &  $ \gamma_{AB} = 0$ \\
 $ \Gamma_{AC} = \{3\}  $ & $\implies$ &  $ \gamma_{AC} = 1$ \\
 $ \Gamma_{BC} = \{4\}  $ & $\implies$ &  $ \gamma_{BC} = 1$ \\
 $ \Gamma_{ABC} = \{1, 2\}  $ & $\implies$ &  $ \gamma_{ABC} = 2$ \\
 $ \Gamma_{\phi} =  \emptyset  $ & $\implies$ &  $ \gamma_{\phi} = 0$ 
 \end{tabular}
 \end{tabular}
 \end{center}
\caption{A partition of the graph union of $\family{C} = \left\{A, B, C\right\}$ with $A = \{1, 2, 3, 5, 6\}, B = \{1,2,4,7,8\}$, and $C =\{1,2,3,4,9\}$ according to its $\Gamma$ sets, together with their sizes $\gamma$.\\
This is also a decomposition of $\Nset{9}$ following Proposition \ref{prop:gamma_partition}.}
\label{fig:three_K5_example}
\end{figure}

Consider again the example of Section \ref{section:intro}, where the collection 
$\family{C} = \{A, B, C\}$ consisting of the three $5$-cliques $A = \{1, 2, 3, 5, 6\}, B = \{1, 2, 4, 7, 8\}$, and $C = \{1, 2, 3, 4, 9\}$ in some graph.
 Figure \ref{fig:three_K5_example} shows the graph union over the cliques of $\family{C}$.

Because various intersections of the cliques in $\family{C}$ are important to identify,  we introduce a separate notation to distinguish those subgraphs, of the graph union over $\family{C}$,  that \textit{uniquely} appear in an intersection of specified cliques in $\family{C}$ \textit{but not} in any of the unspecified cliques.  

\begin{figure}[htbp]
\begin{center}
\begin{tabular}{lc}
 \begin{tabular}{lcl}
  $\Gamma_{A}$ & $ =  $ & $A \intersect \comp{B} \intersect \comp{C} $\\
  $\Gamma_{B}$ & $ = $ & $\comp{A} \intersect B \intersect \comp{C}$ \\
  $\Gamma_{C}$ & $ =  $ & $\comp{A} \intersect \comp{B} \intersect C$  \\
  $\Gamma_{AB}$ & $ =  $ & $A \intersect B \intersect \comp{C}$ \\
  $\Gamma_{AC}$ & $ =  $ & $A \intersect \comp{B} \intersect C $ \\
  $\Gamma_{BC}$ & $ =  $ & $\comp{A}  \intersect B \intersect C$ \\
  $\Gamma_{ABC}$ & $ =  $ & $A \intersect B \intersect C$ 
\\
  $\Gamma_{\phi}$ & $ =  $ & $\comp{A} \intersect \comp{B} \intersect \comp{C}$ 
\end{tabular}
&
\begin{tabular}{c}
\includegraphics[width = 0.3\textwidth]{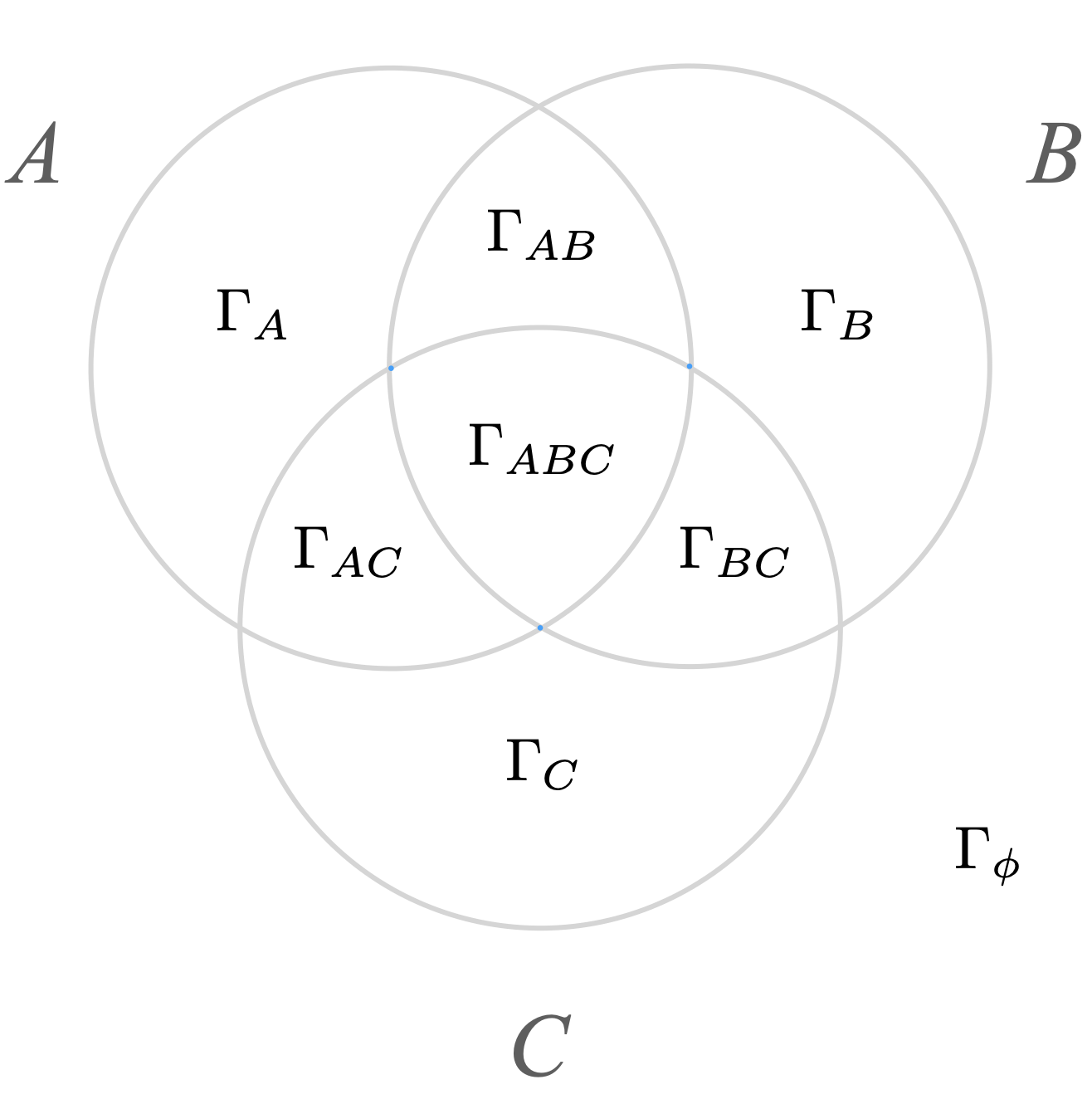}
\end{tabular}
\end{tabular}
\caption{The $\Gamma$ sets for $\family{C} = \{A, B, C\}$.}
\label{fig:gammasets}
\end{center}
\end{figure}

The set of indices is denoted by $\Gamma$ with the specified cliques identified by subscript are shown in Figure \ref{fig:gammasets} 
for a collection of cliques $\family{C} = \{A, B, C\}$.
The $\Gamma$ sets partition its graph union while the indexing on $\Gamma$ partitions the power set of $\mathcal{C}$, which we denote by $\{\emptyset, A, B, C, AB, AC, BC, ABC\}$.

For each cell in $\Gamma$, its cardinality is denoted by $\gamma = \cardinality{\Gamma}$ -- e.g., $\gamma_{AB} = \cardinality{\Gamma_{AB}}$. 
Figure \ref{fig:three_K5_example}
shows the partition of the graph union of $\family{C}$ from Figure \ref{fig:clique_collection} according to its $\Gamma$ sets, as in Figure \ref{fig:clique_collection}. 
The contents of each $\Gamma$ set are easily read off from the graph, as shown. The $\gamma$s are simply the cardinalities of the sets.
%
%
For example, the cell $\Gamma_{AB}$ contains no nodes from the collection because every element common to both $A$ and $B$ is also common to $C$.

The $\Gamma$-sets turn out to have useful properties related to cliques.  
From Figures \ref{fig:gammasets} and \ref{fig:three_K5_example},  note that each original clique $A = \Gamma_{A} \union \Gamma_{AB} \union \Gamma_{AC} \union \Gamma_{ABC}$, $B = \Gamma_{B} \union \Gamma_{AB} \union \Gamma_{BC} \union \Gamma_{ABC}$, and $C = \Gamma_{C} \union \Gamma_{AC} \union \Gamma_{BC} \union \Gamma_{ABC}$, is the union of $\Gamma$-sets whose subscript sets have a common intersection, namely $A$, $B$, or $C$.
Moreover, the size of each clique is simply the sum of the corresponding $\gamma$s.  
Similar results hold for the union of any two $\Gamma$-sets 
$\Gamma_{J_1}$ and $\Gamma_{J_2}$.  If the index sets are such that $J_1\intersect J_2 \ne \emptyset$, then the union $\Gamma_{J_1} \union \Gamma_{J_2}$ forms a clique of size $\gamma_{J_1} + \gamma_{J_2}$; if $J_1\intersect J_2 = \emptyset$, then there is no clique spanning $\Gamma_{J_1}$ and $\Gamma_{J_2}$.

\subsection{An orbit partition}
\label{section:an-orbit-partition}
Consider any $\Gamma$-set in Figure \ref{fig:gammasets} and the node numbers it contains in Figure \ref{fig:three_K5_example}.  
The node numbers within \textit{any} $\Gamma$-set could be permuted without any change in the structure of the graph in Figure \ref{fig:three_K5_example}.  These cells are called \textit{orbits} and the partition an \textit{orbit partition} \cite<e.g., see>[Definition 9.3.4 and Proposition 9.3.5] {rolesLerner05}.
That the $\Gamma$-sets, as defined above, form an orbit partition in general will be proved in Proposition \ref{prop:orbit_partition}.

For any equitable partition (e.g., an orbit partition), $\Gamma = \{\Gamma_1, \ldots, \Gamma_m \}$, of the vertex set of a graph $G$, a directed multi- (or weighted) \textit{quotient} graph can be defined having nodes $\Gamma_i$ and  $b_{ij}$ edges (or edge weights) from  $\Gamma_i$ to $\Gamma_j$ where $b_{ij}$ is the number of neighbours in $\Gamma_j$ of every vertex in $\Gamma_i$ -- called the \textit{quotient} of $G$ modulo $\Gamma$ and denoted $G/\Gamma$ \cite<e.g.,>[Definition 9.3.2]{rolesLerner05}.

\begin{figure}[hbt]
{\scriptsize 
\begin{center}
\begin{tabular}{cr}
\begin{tabular}{c}
\includegraphics[width = 0.2\textwidth]{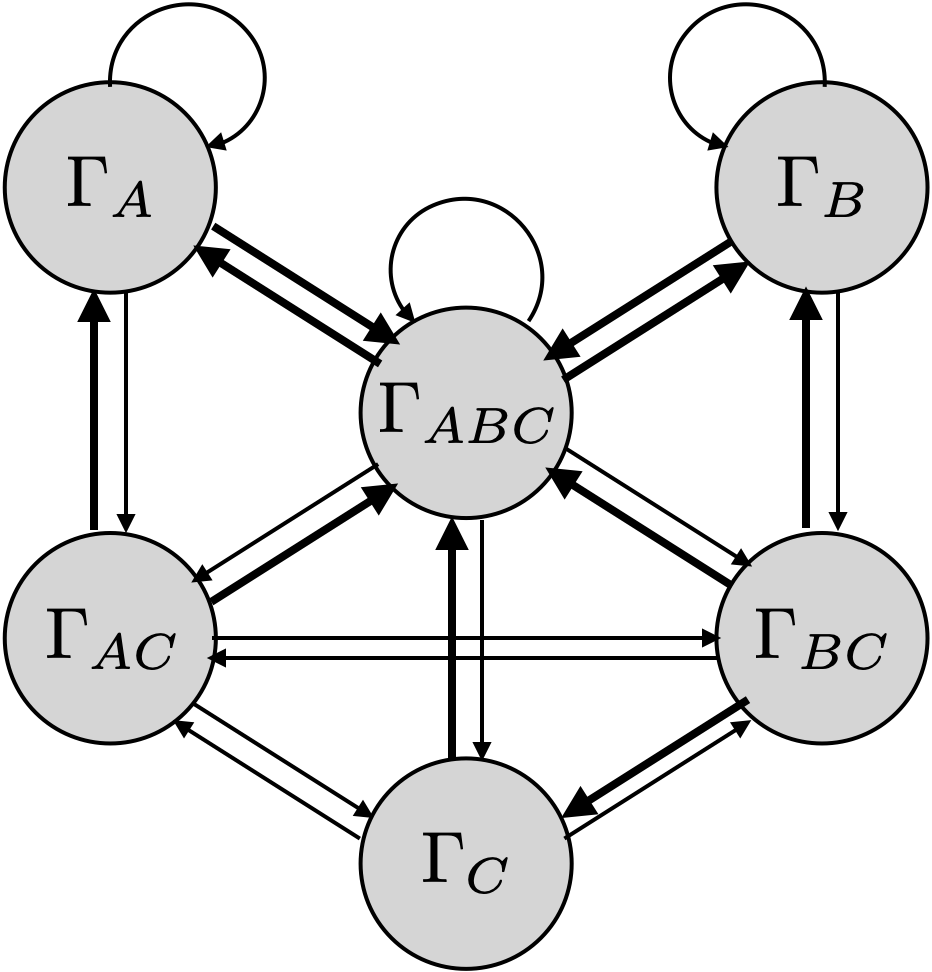}  
\end{tabular}
&
$
\begin{array}{rcl}

orbit & 
~~\Gamma_{A}~~
 ~~\Gamma _{B} 
 ~~~\Gamma_{C} 
 ~~~\Gamma_{AC}  
 ~~~\Gamma_{BC}   
 ~~~\Gamma_{ABC}  
 &  \\
 \begin{array}{l}
\Gamma_{A} \\  
 \Gamma _{B} \\ 
 \Gamma_{C} \\  
  \Gamma_{AC} \\ 
 \Gamma_{BC} \\  
 \Gamma_{ABC} \\ 
 \end{array}
 &
\left[
\begin{array}{cccccc}

   ~~ 1  ~~      &~~   0  ~~&~~   
                                                    0   ~~&~~    1  ~~&~~    0  ~~&~~    2   ~~   \\ 
        0              &      1       &       
                                                    0       &         0      &        1      &         2      \\     
        0              &      0       &       
                                                    0       &         1      &        1      &         2      \\     
        2              &      0       &       
                                                    1       &         0      &        1      &         2      \\     
        0              &      2       &       
                                                    1       &         1      &        0      &         2      \\     
        2              &      2       &       
                                                    1       &         1      &        1      &         1   \\     

\end{array}
\right]
& 
 = \m{B} \end{array}
$

\end{tabular}
\end{center}
}
\caption{The quotient graph of the graph union of $\family{C}$ modulo $\Gamma$ and its edge weight matrix $\m{B} = [b_{ij}]$. Edges are shown with width proportional to their weight in $\m{B}$.}
\label{fig:quotient_graph}
\end{figure}

For the graph union of Figure \ref{fig:three_K5_example}, the partition $\Gamma = \{\Gamma_A, \Gamma_B, \Gamma_C, \Gamma_{AC}, \Gamma_{BC}, \Gamma_{ABC}\}$ produces the quotient graph and matrix $\m{B} = [b_{ij}]$ shown in Figure \ref{fig:quotient_graph}.
This graph can be thought of as a compression of the original graph union.  As such, some information will be lost, but much remains. 
Its (weighted) adjacency matrix and graph are enough to determine several properties of the graph union \cite<e.g., see >{godsil1993algebraic}, including the path distances between nodes, the graph diameter, and a partial spectral decomposition -- the characteristic roots of $\m{B}$ are a subset of those of the adjacency matrix $\m{A}$ of the graph union.

\subsection{Equivalent graphs}
\label{section:equivalentgraphs}
The orbit partition, $\Gamma$, has particular features.
For example,  
if \textit{any} node $u$ in $\Gamma_{J_1}$ connects to $k$ nodes in  $\Gamma_{J_2}$, then \textit{every} node in  $\Gamma_{J_1}$ connects to the \textit{same} $k$ nodes in  $\Gamma_{J_2}$, and vice versa.
And, since permuting node numbers in any $\Gamma$-set does not change the graph, if $u \sim v$ for any $u \in \Gamma_{J_1}$ and any $v  \in \Gamma_{J_2}$, then \textit{all} nodes in $\Gamma_{J_1}$ connect to \textit{all} nodes in $\Gamma_{J_2}$. 
It follows, then, that the nodes of $\Gamma_{J_1} \union \Gamma_{J_2}$ form a clique of size $\gamma_{J_1} + \gamma_{J_2}$ (since each of $\Gamma_{J_1}$ and $\Gamma_{J_2}$ also form cliques).

This suggests that the orbit partition, given by the $\Gamma$-sets,  provides a structure to identify sets of equivalent subgraphs which may, or may not, form a clique. 
If we choose an ordering of the orbits, say $(\Gamma_{A}, \Gamma_{B}, \Gamma_{AB}, \Gamma_{C}, \Gamma_{AB}, , \Gamma_{AC}, \Gamma_{ABC})$, then a unique tuple of the counts of nodes from each orbit identifies a set of subgraphs which are isomorphic to one another (under node permutation within each orbit).   
For example, both $\{1, 3, 6\}$ and $\{2, 3, 6\}$ share the tuple 
(1, 0, 0, 0, 1, 0, 1), but $\{1, 2, 3\}$ with tuple (0, 0, 0, 0, 1, 0, 2) is a unique subgraph (under permutation within orbits).
Each of these forms a $3$-clique.  The size of the subgraph is the sum of the tuple elements and the number of subgraphs the tuple represents is the product of the size of the orbit choose that element of the tuple  (e.g., here $\binom{2}{1} \times \binom{1}{1} \times \binom{2}{1} = 4$ in total, the remaining two being $\{1, 3, 5\}$ and $\{2, 3, 5\}$).

The index sets associated with each non-zero tuple element determine whether subgraphs produced by the tuple are also a clique.
For example, the tuple (1, 0, 0, 0, 1, 0, 1) takes nodes from $\Gamma_{A}$,  $\Gamma_{AC}$, and $\Gamma_{ABC}$ whose index sets are $\{A\}$, $\{A, C\}$, and $\{A, B, C\}$.   The intersection of these index sets is non-null, and every subgraph induced by this tuple is a clique of size equal the sum of its elements.  In contrast, the index sets $\{B\}$ and $\{C\}$ corresponding to the tuple (0, 2, 0, 1, 0, 0, 0) have null intersection and this tuple's induced graph does not form a clique.  A tuple induces a clique if, and only if, the intersection of its index sets is non-null -- this is formally established by Proposition \ref{prop:adjacency_intersection}. The tuple associated with a clique we call its \textit{signature}, and cliques having the same signature are of the same \textit{type}.

\subsection{Maximal cliques}
\label{section:maximal-example}
Consider the problem of finding all maximal cliques which contain some specific clique. 
For example, from Figure \ref{fig:three_K5_example}, find all maximal cliques which contain the $2$-clique $\{1, 2\}$.  These are
\begin{enumerate}[label = (\roman*)]
	\item $M_1 = \{1, 2, 3, 5, 6\} = \Gamma_{A} \union \Gamma_{AC} \union \Gamma_{ABC}$,
	\item $M_2 = \{1, 2, 4, 7, 8\}  = \Gamma_{B} \union \Gamma_{BC} \union \Gamma_{ABC}$, and
	\item $M_3 = \{1, 2, 3, 4, 9\} = \Gamma_{C} \union \Gamma_{AC} \union \Gamma_{BC} \union \Gamma_{ABC}$.
\end{enumerate}
The maximal cliques help enumerate the total number of cliques which contain a specified clique by identifying the nodes which can be added to expand that clique. 
In the case of $\{1, 2\}$, $M_1$ provides three additional nodes (viz., 3, 5, and 6) and so $2^3 - 1 = 7$ larger cliques containing $\{1, 2\}$. The same holds for $M_2$ and $M_3$, but care must be taken for double counting.  The total number of cliques containing $\{1, 2\}$ (including itself) is expressed in terms of its maximal cliques as
\begin{align*}
&\sum_{i=1}^{3}\left(2^{|M_i|- |\{1,2\}|}-1\right) 
  ~-  \sum_{\{i, j\} \subset \{1,2,3\}} \left(2^{|M_i\intersect M_j|- |\{1,2\}|}-1 \right)
  ~+ ~\left(2^{|M_1\intersect M_2 \intersect M_3| - |\{1,2\}|}-1\right) 
  + 1\\ 
 &\\
&~~~~~~~~~~= (7 + 7 + 7)  -  \left((2^0 - 1) +(2^1 - 1) + (2^1 - 1)   \right)+ (1 -1) +1 = 20,
\end{align*}
where the last summand 1 corresponds to the edge $\{1,2\}$ on its own. 
A general expression for this count is given in Proposition \ref{prop:cliques_containing_H}.

We might call the union $M_1 \union M_2 \union M_3$, of its maximal cliques, the  \textit{clique extension} of $\{1, 2\}$ within the cover, or simply the \textit{clique extent} of $\{1, 2\}$.  In this case, the extent of $\{1, 2\}$ is the entire cover but this is not generally the case (e.g., the extent of $\{5, 6\}$ is simply the set $A = \{1, 2, 3, 5, 6\}$).  
More generally, if the intersection of all cliques in the collection is non-null, then the clique extension of any node (or clique) in that intersection will generate the entire cover.  
In a social network context, for example, such individuals (or cliques) might be deemed to be highly influential in the entire cover -- wherever they are located in the cover, those having larger clique extents might be regarded as more influential than those having smaller ones.

\subsection{Intersecting families}
\label{section:intersectingfamilies-example}
Reading off the set of subscripts from the $\Gamma$-sets defining each of the maximal cliques, $M_1, M_2, M_3$, respectively, gives the sets:
\begin{enumerate}[label = (\roman*)]
	\item $\family{F}_1 = \{\{A\}, \{A, C\}, \{A, B, C\} \}$,
	\item $\family{F}_2 = \{\{B\}, \{B, C\},  \{A, B, C\} \}$,
	\item $\family{F}_3 = \{\{C\}, \{A, C\}, \{B, C\}, \{A, B, C\} \}$.
\end{enumerate}
Each of these sets, $\family{F}_i$, is called an \textit{intersecting family}  \cite<e.g., see>{meyerowitz1995maximal}, meaning that each set  $\family{F}_i$ is a subset of the power set of $\{A,B,C\}$ and that its elements have non-null pairwise intersection.  
When the context is clear,  the notation for an intersecting family will be simplified,  from a set of sets, to a set of the subscripts identifying the corresponding $\Gamma$-sets --  so, the contents of $\family{F}_1$ can be simplified to $\{A,  AC,  ABC \}$.
Being cliques, each of the above families share the additional property that they have non-null intersection over all of their elements, not just pairwise. 

When an intersecting family $\family{F}$ is not a proper subset of any other intersecting family, it is called a \textit{maximally intersecting family} \cite{meyerowitz1995maximal}. The family $\family{F}_3 = \{C, AC, BC, ABC\}$ is a maximally intersecting family and corresponds to the maximal clique $c_3$.  The families $\family{F}_1$ and $\family{F}_2$ are not; though, since $\Gamma_{AB} = \emptyset$, adding $AB$ to each will make them maximally intersecting as well (i.e., equivalently, $M_1 = \Gamma_{A} \union \Gamma_{AB} \union \Gamma_{AC} \union \Gamma_{ABC}$ and $M_2 = \Gamma_{B} \union \Gamma_{AB} \union \Gamma_{BC} \union \Gamma_{ABC}$).
Note that maximal intersecting families are not necessarily isomorphic to maximal cliques.
For example, the only other maximal intersecting family here is $\family{F}_4 = \{AB, AC, BC, ABC\}$ corresponds to the $4$-clique $\{1, 2, 3, 4\}$ which is not maximal.
Necessary and sufficient conditions for an intersecting family to determine a maximal clique are given in Theorem \ref{thm:maximal_clique} provides the necessary and sufficient conditions for a clique to be maximal.
Intersecting families have interesting structure 
\cite<e.g., see>{meyerowitz1995maximal}.

We call a collection of sets $\mathcal{F}$  \textit{path intersecting} if, for any $A, B \in\mathcal{F}$ there exists a sequence sets  $J_1, J_2,\ldots, J_\ell$ in $\mathcal{F}$ from $A = J_1$ to $B = J_\ell$ having  that $J_{j}\cap J_{j+1} \neq \emptyset$ for all $j=1,2,\ldots, \ell - 1$. 
The set collection $\{ \{A\}, \{B\}, \{C\}, \{A, B\}, \{A, C\}\}$ is path intersecting, but is not an intersecting family.  

In Section \ref{section:generalquotient}, the sequence of index sets of nodes along any path in the quotient graph is shown to be path intersecting and the set of index sets from any clique on the same to be an intersecting family.

%

\section{The general approach} 
\label{section:partitional-framework}
By example, a number of results were illustrated in Section \ref{section:partition-example} relating the properties of a particular partition of the vertex set of the graph union of a collection of cliques to the cliques within the union. In this section, we show more generally that this kind of  partition is a link between cliques in the union of a clique collection and certain intersecting families of sets. This link allows for a nuanced enumeration of several clique counting problems on these graphs, including total number of cliques, maximal cliques, maximum cliques and cliques containing any specific subset of interest.  Finally, as with the example of Section \ref{section:partition-example}, we establish that this kind partition is an orbit partition, hence capturing salient features of the original graph.  

The example clique collection of Section \ref{section:partition-example} had three maximal $5$-cliques as its elements  -- while possibly desirable, this is not necessary. 
In this section, general results for cliques in the graph union of \textit{any} collection of cliques are derived (i.e., each of any size, including possibly as a single edge).
We begin with the general construction of a vertex partition of the graph union which permits a deeper examination of all cliques through intersecting families derived from that partition. Maximal intersecting families will be shown to correspond to the largest cliques obtainable from particular sub-collections of cliques. 

\subsection{The partition}
\label{section:partition-general}
The general construction of the partition, and its cells,  are defined in Proposition \ref{prop:gamma_partition}.
Here, for any set $J \subseteq \Nset{m}$, its set complement is with respect to $\Nset{m}$ and is denoted as $\comp{J} := \Nset{m} \setminus J$.
\begin{proposition}
\label{prop:gamma_partition}

For any $m\geq 1$, given a sequence $(A_i)_{i=1}^{m}$ of subsets of $\Nset{n} =  \union_{i=1}^m A_i$, the family of sets given by

$$\Gamma := \left\{\bigcap_{i\in J}A_i\setminus \left(\bigcup_{i  \in \comp{J}}A_i \right): J\subseteq \Nset{m} \right\} :=  \left\{\Gamma_J : J\subseteq \Nset{m} \right\}$$
is a partition of  $\Nset{n}$. 

\noindent
Moreover, for any $i \in \Nset{m}$, 
$$A_i = \bigcup_{J\subseteq \Nset{m} :~ i \in J } \Gamma_J.$$
\end{proposition}
\begin{proof}
First, we show that 
\[
 \bigcup_{J\subseteq \Nset{m}} \Gamma_J = \bigcup_{J\subseteq \Nset{m}}
					\left[ \bigcap_{i\in J}A_i\setminus 
							\left(\bigcup_{i \in \comp{J}} A_i \right)
							\right] = \Nset{n}.
							\]
\noindent
For every $J\subseteq \Nset{m}$, 
\[\Gamma_J = \left[ \bigcap_{i\in J}A_i\setminus 
							\left(\bigcup_{i \in \comp{J}} A_i \right)
							\right] \subseteq\Nset{n},
							\]
as each $A_i \subseteq \Nset{n}$. To see the reverse inclusion, fix \textit{any} choice $x \in \union_{i=1}^m A_i = \Nset{n}$ and
let $J_x := \{ i: x \in A_i\}\subseteq \Nset{m}$ denote the set of all indices $i$ with $x\in A_i$, and its complement in $\Nset{m}$ as $\comp{J_x} = (\Nset{m}\setminus J_x)$.
Now $x \in \Nset{n}$ appears in at least one $A_i$, since $\union_{i=1}^m A_i = \Nset{n}$,
so it follows that $x \in \intersect_{i \in J_x} A_i$  and $x \not\in \union_{i \in \comp{J_x}} ~A_i$. 
Thus,
\[
x \in \left[ \bigcap_{i\in J_x}A_i\setminus \left(\bigcup_{i \in \comp{J_x}}A_i \right)\right]
= \Gamma_{J_x}
\]
for \textit{any} $x \in \Nset{n}$, and hence 
\[
\Nset{n}
= \bigcup_{x \in \Nset{n}} \Gamma_{J_x} =  \bigcup_{J\subseteq \Nset{m}}
                  \Gamma_J.
							\]				
It remains only to show that the intersection of any two distinct non-null members of $\Gamma$ is empty -- the proof is by contradiction.
Let $J, H\subseteq \Nset{m}$  be distinct, respectively producing  
\[
\Gamma_J = \left[ \bigcap_{i\in J}A_i\setminus 
		\left(\bigcup_{i\not\in J}A_i \right)\right] ~~\text{and}~~
				\Gamma_H = \left[ \bigcap_{i\in H}A_i
					\setminus \left(\bigcup_{i\not\in H}A_i \right)
					\right]
\]
as members in $\Gamma$. 
Suppose $x \in \Gamma_J \intersect \Gamma_H \ne \emptyset$, then $x \in\Gamma_J \imply x \in A_i ~ \forall i \in J$ and $x \in\Gamma_H \imply x \in A_i ~ \forall i \in H$. 
Since $J$ and $H$ are distinct, there exists some $k\in J\setminus H$ for which $x \in \Gamma_J $ appears in $A_k$.  Now $k \not \in H$ means $k \in \comp{H}$ and hence $A_k$ appears in the union $\union_{i\not\in H}A_i$ being removed from  $\intersect_{i\in H} A_i$ in the definition of $\Gamma_H$.  Therefore $x \not \in \Gamma_H$ and, so, $x \not \in \Gamma_J \intersect  \Gamma_H$, a contradiction.  
It follows that $\Gamma_J $ and $\Gamma_H$ are disjoint, whenever $J \ne H$ and hence that the sets of $\Gamma$ form a partition of their union, $\Nset{n}$.

Finally, for any $i \in \Nset{m}$, it remains only to show that the original sets $A_i$ are the union of those $\Gamma$-sets, $\Gamma_J$, whose index set $J$ contains $i$.  That is, 
$$A_i = \bigcup_{J\subseteq \Nset{m} :~ i \in J } \Gamma_J. $$
If $i \in J$, then $\Gamma_J = \left[ \bigcap_{j\in J}
        A_j\setminus 
	\bigcup_{j \in \comp{J}}A_j \right]$ intersects $A_i$, and hence $\Gamma_J \subseteq A_i$ whenever $i \in J$.  It follows, then, that
\[
 \bigcup_{J\subseteq \Nset{m} :~ i \in J } \Gamma_J \subseteq A_i.
 \]
Conversely, for every $x \in A_i$, then $i \in J_x$ and 
\[
x \in 
   \left[ \bigcap_{j\in J_x}A_j \setminus \bigcup_{j \in \comp{J_x}}A_j \right] 
  = \Gamma_{J_x} \subseteq \bigcup_{J\subseteq \Nset{m} :~ i \in J } \Gamma_J.
\]
So $A_i \subseteq \bigcup_{J\subseteq \Nset{m} :~ i \in J } \Gamma_J \subseteq A_i$, and
it follows that $A_i = \bigcup_{J\subseteq \Nset{m} :~ i \in J } \Gamma_J$.
\end{proof}
We will call a partition produced as in Proposition \ref{prop:gamma_partition}, a \textit{$\Gamma$-partition} and note that it will be peculiar to the sets $A_i$ from which it is constructed.

\subsubsection{Applied to a clique collection} 
\label{section:cliquepartition}


For a collection of cliques $\mathcal{C} = \left\{ c_1,  \ldots, c_m \right\}$, defined by index sets $c_j \subset \Nset{n}$, with graph union $\bigcup_{j = 1}^m c_j = \Nset{n}$,  Proposition \ref{prop:gamma_partition} provides a general means to find $\Gamma$-sets, namely as $(\Gamma_J)_{J\subseteq \Nset{m}}$ with 
$$\Gamma_J = \left(\bigcap_{j\in J}c_j\right)\intersect \left(\bigcap_{j\not\in J}\comp{c_j}\right)$$
where complement is with respect to $\Nset{n}$.
That is, each cell $\Gamma_J$ is the set of vertices common to all $c_j$ for all $j\in J$ and absent from every $c_j$ for which $j\not\in J$. 
Again, the cardinality of $\Gamma_J$ is denoted as $\gamma_J = |\Gamma_J|$.  

The
$\Gamma$-partition provides
an equivalence relation on nodes 
 $u, v \in \Nset{n}$ via the indices of those cliques which contain $u$ or $v$
 -- namely, $J_u = \{j \in \Nset{m}: u \in c_j\}$ and $J_v = \{j \in \Nset{m}: v \in c_j\}$.
The nodes $u$ and $v$ are equivalent, $u\equiv v$, if, and only if, $J_u = J_v$; that is, $u$ and $v$ are in the same $\Gamma$-set.

The $\Gamma$-partition can also be used directly to infer some properties of the graph union.
For example, as in Section \ref{section:equivalentgraphs}, the adjacency of nodes in the graph union is related to the intersection of those $\Gamma$-sets which contain them: 
\begin{proposition}
\label{prop:adjacency_intersection}
Let $u$ and $v$ be two nodes in the graph union, $\bigcup_{j=1}^{m}c_j$,  of the clique collection $\mathcal{C} = \left\{ c_1, c_2, \ldots, c_m \right\}$. If $u \in \Gamma_{J_u}$ and $v\in \Gamma_{J_v}$, then $u\sim v$ if, and only if, $J_u \intersect J_v \ne \emptyset$.
\end{proposition}
\begin{proof}
	We note that $u \sim v$ if, and only if, for some $j \in \Nset{m},$ $u\in c_j$ and $v\in c_j$, which is equivalent to $J_u \intersect J_v \neq \emptyset$.
\end{proof}
\noindent
It follows, for example, that $u \sim v$ for every pair of nodes $u, v \in \Gamma_J$ (for any $J \subseteq \Nset{m}$).
Moreover, the cardinalities, $\gamma_J$, determine the degree of every vertex in $\Gamma_J$. 
Proposition \ref{prop:degree_seq_in_cliques} establishes that all nodes in a cell of $\Gamma$ have the same degree. 
\begin{proposition}
\label{prop:degree_seq_in_cliques}
For a non-null set $J \subset \Nset{m}$, every vertex in $\Gamma_J$ has degree $d_J$ where
$$d_J= \sum_{I \subseteq \Nset{m} ~:~ I \intersect J \ne \emptyset}\gamma_I   ~~~~~ - 1.$$
\end{proposition}
\begin{proof}
If $u\in \Gamma_J$, then $u\sim v$ if, and only if, 	
$$v \in \bigcup_{I \subseteq \Nset{m} ~:~ I \intersect J \ne \emptyset} \Gamma_I = \bigcup_{j\in J}c_j,$$
with $v\neq u$. Therefore, the degree of $u$ is
\begin{align*}
  deg(u) &= \left|\bigcup_{j\in J}c_j \right| - 1\\
  &= \left|\bigcup_{j\in J}\left(\bigcup_{I: j\in I}\Gamma_I \right)\right| - 1\\
  &= \sum_{I: j \in I, \text{ for some } j \in J} \gamma_I - 1.
\end{align*}
\end{proof}

Note that different clique collections having the same graph-union produce different $\Gamma$-partitions, these being peculiar to the particular cliques in the collection.  
The cliques of the collection in Section \ref{section:partition-example}, for example, were all of size 5; had they all been of size 3 the same graph union of (now many more) cliques in the collection would be the same but the resulting $\Gamma$-sets would be different. 

The special case that the collection consists of exactly $m$ cliques of size $r$, as in Section \ref{section:partition-example}, can also be determined from the cardinalities, $\gamma_J$.
For $\Gamma$ to have been formed from a collection of $m$ distinct $r-$cliques, the following must hold:
\begin{align*}
	\sum_{J\subseteq\Nset{m}}\sv{\gamma}_J &= n, 
	    \tag*{\ldots for the graph union to have $n$ nodes}\\
	\sum_{J\subseteq\Nset{m} ~:~ j\in J} \sv{\gamma}_J &= r
	    \tag*{\ldots for each $c_j$ to have $r$ nodes}\\
	\sum_{J\subseteq\Nset{m} ~:~  \{j,k\} \subseteq J} \sv{\gamma}_J &< r
	    \tag*{\ldots to ensure distinct cliques: when $j \ne k$,  $c_j\neq c_k$.}
\end{align*}
Of these, only the last may not be self-evident; it follows from:
\begin{proposition}
Let $\mathcal{C} = \mCliqueSet{c}$ be a collection of $r-$cliques and fix $I \subseteq \Nset{m}$. Then $\{c_i: i\in I\}$ consists of a single clique if, and only if, \[\sum_{J: I\subseteq J} \gamma_J = r.\] 	
\end{proposition}
\begin{proof}
By Proposition \ref{prop:gamma_partition}
\[ \bigcap_{i \in I} c_i = \bigcap_{i\in I} \bigcup_{J\subseteq \Nset{m} : i \in J}\Gamma_J = \bigcup_{J:I \subseteq J}\Gamma_J \]
So, $|\cap_{i\in I} c_i |= |\bigcup_{J:I \subseteq J}\Gamma_J| = \sum_{J: I\subseteq J} \gamma_J.$ Since all $c_i$ are $r-$sets, their intersection is an $r-$set if, and only if, they are all equal. 
\end{proof}

\subsection{The general $\Gamma$-quotient graph}
\label{section:generalquotient}

In light of Proposition \ref{prop:adjacency_intersection}, $\Gamma$ is an  \emph{equitable partition} \cite<e.g., see>{godsil2001algebraic, rolesLerner05} -- the number of neighbours in $\Gamma_H$ of vertex $u \in \Gamma_J$ depends only on the choice of $H$ and $J$. In fact, $\Gamma$ is an orbit partition induced by a group of automorphisms of $H$.

\begin{proposition}
	\label{prop:orbit_partition}
	The partition $(\Gamma)_{\emptyset \neq J \subseteq \Nset{m}}$ is an orbit partition. 
\end{proposition}
\begin{proof}
For a nonempty $J\subseteq \Nset{m}$, let $\pi_J$ be any permutation of the elements of $\Gamma_J$. Let $\pi:V\to V$ be the extension of the $\pi_J$ to $V$. It immediately follows that the orbits of $\pi$ are the cells of $\Gamma$ and, by Proposition $\ref{prop:adjacency_intersection}$, that $\pi$ is an automorphism of $V$.
\end{proof}
Each $\Gamma_J$ cell has an $n \times 1$ characteristic  vector $\ve{c}_J$ having value $1$ in row $i$ if vertex $i$ is in $\Gamma_J$, and $0$ otherwise, so that $\tr{\ve{c}_J} \ve{c}_J = \gamma_J$.  The characteristic matrix $\m{C}$ is formed with columns $\ve{c}_J$ placed in order of the $\Gamma_J$s of the partition $\Gamma$.
If $\m{A}$ is the adjacency matrix of the graph union, $G$, the matrix $\m{B} = \inv{(\tr{\m{C}}\m{C})}  \tr{\m{C}} \m{A} \m{C}$ determines the structure of the quotient graph of $G$ modulo $\Gamma$ \cite<e.g., see>[Lemma 9.3.1, p. 196]{godsil1993algebraic}.

\subsection{Type equivalent graphs}
\label{section:types}
For the example of Figure \ref{fig:three_K5_example}, Section \ref{section:equivalentgraphs}, introduced the \textit{type} of a subgraph $H$ of the graph union of $\mathcal{C} = \mCliqueSet{c}$ associated with its $\Gamma$-partition and identified by a \textit{signature}, namely,  the tuple of the counts of nodes from $H$ appearing in each cell of the partition. 
In this section, these ideas are formalized to provide a more nuanced sense of equivalent graphs in the context of a $\Gamma$-partition of the graph union  $G = \union_{j = 1}^m c_j$.

For subgraph $H$ of $G$, the \textit{signature} of $H$ defined by the $\Gamma$-partition of $\mathcal{C}$ is the \textit{function} $f_H: \powerset{\Nset{m}} \to \NaturalsZero$ defined as $f_H(J) = |H \cap \Gamma_J|$ for all $J \subseteq \Nset{m}$.  Note that this is defined for any subgraph $H$, not necessarily only cliques $H$.  Two subgraphs $H_1$ and $H_2$ are said to be of the same \textit{type}, or to be \textit{type-isomorphic}, if, and only if, they have identical signatures (i.e., $f_{H_1} = f_{H_2}$).  Finally, the \textit{support} of $H$ (or of $f_H$) is the set of all subsets $J$ of $\Nset{m}$ for which $f_H(J) > 0$; we write the support as $\support{H} = \{J : J \subseteq \Nset{m} \text{ and } f_H(J)  > 0 \}$, or as $\support{f_H}$ when emphasizing the signature.  Note also that all of these are predicated on the particular clique collection $\mathcal{C}$ and its associated $\Gamma$-partition.

For example, consider the  clique collection of Figure \ref{fig:three_K5_example} and the subgraphs $H_1 = \{1, 2, 3, 4\}$, $H_2 = \{1, 2, 3, 5\}$, $H_3 = \{1, 2, 3, 6\}$, and $H_4 = \{1, 2, 3, 5, 6\}$.  The first three are graph isomorphic to each other and the complete graph, $K_4$ while $H_4$ is isomorphic to $K_5$.
In contrast only $H_2$ and $H_3$ are type isomorphic; $H_1$ has a different signature (and support), while $H_4$ shares the same support as $H_2$ and $H_3$ but is of a different type.

Because it differs from the usual graph equivalence, the notion of type could be of interest whenever the node labels, or the cliques defining the collection, carry additional meaning.

\subsubsection{$\Gamma$-signatures}
\label{section:signatures}
This section develops a number of counting results obtained types of subgraphs (as defined by signature) from any specific clique collection.

The number of different types of induced subgraphs is easily captured by the cell sizes of the partition:

\begin{proposition}
\label{prop:type_isomorphism_classes}
The number of 
distinct signatures for the
$\Gamma$-partition of 
a collection of $m$ cliques is 
\[ \prod_{J\in \powerset{\Nset{m}}} (\gamma_J+1). \]
\end{proposition}
\begin{proof}
A function $f:\powerset{\Nset{m}} \to \NaturalsZero$ is a signature if, and only if, $|f(J)| \leq \gamma_J$. Thus, there are $\gamma_J+1$ choices for every $J \in \powerset{\Nset{m}}$. 
\end{proof}

\begin{proposition}
\label{prop:num_signatures_same_support}
For any signature $f_H$, the number of signatures having the same support, $\support{H}$, is 
\[
\prod_{J \in \support{H}} \gamma_J.
 \]
\end{proposition}
\begin{proof}
For signatures $f_{H_1}$ and $f_{H_2}$ to have the same support, they must have the same $\Gamma$-cells, $\Gamma_J$ for $J \in \support{H_1} = \support{H_2}$, and each signature can have values $1, \ldots, \gamma_J$ for  the $J$th cell.  The total possible is therefore  $\prod_{J \in \support{H}} \gamma_J$.
\end{proof}


\begin{proposition}
\label{prop:num_type_isomorphism_classes}
Let $f: \powerset{\Nset{m}} \to \NaturalsZero$. The number of induced subgraphs having signature $f$ in the graph union of the clique collection $\{c_1,\ldots c_m\}$ is 
\[ 
\prod_{J \in \support{f}} \binom{\gamma_J}{f(J)}.
\]
\end{proposition}
\begin{proof}
The signature is invariant to the choice of nodes within each $\Gamma$-cell -- provided the same number of nodes from each cell is chosen, the signature is the same. Each cell has $\gamma_j$ nodes giving
\[ \prod_{J \in \support{f}}\binom{\gamma_J}{f(J)}\]
choices for type-isomorphic induced subgraphs.
\end{proof}

\subsubsection{Connected subgraphs}
\label{section:connected-subgraphs}
The $\Gamma$-signature of an induced graph also tells whether it is \textit{connected}. This is captured by the notion of a \textit{path-intersecting} collection of sets defined in Section \ref{section:intersectingfamilies-example}.


\begin{proposition}
\label{prop:num_signatures_same_support}
A subgraph $H$ of the graph union over a clique collection 
is connected if, and only if, its support is path-intersecting.
\end{proposition}
\begin{proof}
Since, $f_H$ is defined by the $\Gamma$-partition of $\mathcal{C}$, every node must appear in exactly one set $J$ of $\support{H}$.
Moreover, any pair of nodes $u, v \in H$ appearing in the same set $J \in \support{H}$ are connected by construction of the partition.  So, we need only consider nodes $u$ and $v$ which lie in different sets of the support.  

Suppose $\support{H}$ is path-intersecting.  Then for any pair of nodes $u, v \in H$, which appear in different subsets $J_u, J_v \in \support{H}$, a sequence of sets $J_{w_1}, J_{w_2}, \ldots, J_{w_\ell}$ can be found in $\support{H}$ such that $J_u = J_{w_1}$, $J_v = J_{w_\ell}$, and $J_{w_i} \intersect J_{w_{i+1}} \neq \emptyset$ for all $i = 1, \ldots, (\ell -1)$.  From Proposition  \ref{prop:adjacency_intersection} $w_i \sim w_{i+1}$ for all $i = 1, \ldots, (\ell -1)$, $u = w_1 \rightarrow w_2 \rightarrow \cdots \rightarrow w_\ell = v$, is a path from $u$ to $v$ in $H$, and so the subgraph $H$ is connected.

Conversely, suppose $H$ is connected.  Every pair of nodes $u, v$ appearing in separate sets $J_u$ and $J_v$ of $\support{H}$ have a path connecting them in $H$.  By the construction of $\Gamma$, this path can be chosen to be
$u = w_1 \rightarrow w_2 \rightarrow \cdots \rightarrow w_\ell = v$ such that each $w_i$ comes from a different $J_i$ in $\support{H}$. Again, by  Proposition  \ref{prop:adjacency_intersection},  $w_i \sim w_{i+1}$ implies $J_i \intersect J_{i+1} \neq \emptyset$, and hence that $\{J_1, \ldots, J_\ell\}$ is path-intersecting.  This holds for any $u, v \in H$ and hence any $J_u, J_v \in \support{H}$, implying that it holds for the whole of $\support{H}$.  It follows that $\support{H}$ is path-intersecting. 
\end{proof}

\begin{proposition}
\label{prop:type_isomorphism_classes}
Let $\mathcal{I}_P$ be the set of all path-intersecting collections of non-empty cells from the $\Gamma$-partition of a clique collection $\mathcal{C}$.
The number of distinct signatures that induce a connected subgraph in the graph union over $\mathcal{C}$  
is
\[ \sum_{\mathcal{F}\in \mathcal{I}_P} \prod_{J\in \mathcal{F}} \gamma_J. \]
\end{proposition}
\begin{proof}
Proposition \ref{prop:num_signatures_same_support} states that for a 
subgraph $H$ to be connected, its support must be path-intersecting; Proposition \ref{prop:num_signatures_same_support} determines the number of distinct signatures having the same support.
Together they give the result.
%
%
\end{proof}
\noindent
It follows that the number of induced \textit{disconnected} subgraphs is
\[ \prod_{J\in \powerset{\Nset{m}}} (\gamma_J+1) - 
	\sum_{\mathcal{F}\in \mathcal{I}_P} \prod_{J\in \mathcal{F}} \gamma_J
\]
where $\mathcal{I}_P$ denotes the set of all path-intersecting collections of non-empty cells from $\Gamma$.

\begin{proposition}
\label{prop:gen_fn_induced_connected}
Let $\mathcal{I}_P$ be the set of all path-intersecting collections of non-empty cells from the $\Gamma$-partition of a clique collection $\mathcal{C}$.
The number of induced connected subgraphs of size $k$ in the graph union over $\mathcal{C}$  
is the $k$-th coefficient of the generating series
\[ \sum_{\mathcal{F}\in \mathcal{I}_P} \prod_{J\in \mathcal{F}} \left[(1+x)^{\gamma_J}-1\right]. \]
\end{proposition}
\begin{proof}
By Proposition \ref{prop:num_signatures_same_support}, every induced connected subgraph $H$ is contained in some path-intersecting family. In fact, there exists a unique smallest path-intersecting family $\mathcal{F}_H := \support{H}$ containing it. Clearly, the contribution of $H$ to the generating function
 \[\sum_{H'} x^{|V(H')|}\]
 is $x^k$, where $|V(H)|= k$, and the sum is over all $H'$ induced connected subgraphs whose support is $\mathcal{F}_H$. 
 
 Conversely, given a path-intersecting family $\mathcal{F}$, the induced connected subgraphs whose support is $\mathcal{F}$ are constructed uniquely by choosing $\alpha_J\geq 1$ nodes from $\Gamma_J$ for every $J\in \mathcal{F}$. The generating series corresponding to this is 
 \[\prod_{J\in \mathcal{F}} \left[(1+x)^{\gamma_J}-1\right].  \]
\end{proof}
\subsubsection{$\Gamma$-support and cliques}
\label{section:supportandcliques}
The support of a subgraph $H$ provides information on whether $H$ is a clique and whether it is maximal.
\begin{proposition}
\label{prop:cliques_to_intersecting_families}
For any clique collection $\family{C} =\mCliqueSet{c}$,  the subgraph induced by $H$ on the graph union $\bigcup_{j =1}^{m} c_m$, is a clique, if, and only if, 
is support, $\support{H} = \{J: J \subseteq{\Nset{m}} \text{ and } \Gamma_J \intersect H \ne \emptyset \}$ is an intersecting family.
\end{proposition}
\begin{proof}
Suppose the induced graph on $H$ is a clique. 
Fix two distinct sets $J_1, J_2 \in \support{H}$. 
Let $u_1 \in \Gamma_{J_1}\intersect H$ and $u_2 \in \Gamma_{J_2}\intersect H$. 
Since $u_1 \sim u_2$, it must be that $u_1, u_2 \in c_j$ for some $j\in \Nset{m}$. 
Therefore, it follows that $j \in J_1$ and $j\in J_2$, by the definition of the partition $(\Gamma_J)_{J\subseteq \Nset{m}}$. Thus, $|J_1 \intersect J_2| \geq 1$ and 
$\support{H} $ 
is an intersecting family.

On the other hand, suppose that $\support{H}$
is an intersecting family. Fix $u, v \in H$ and suppose that $u \in \Gamma_{J_u}$ and $v\in \Gamma_{J_v}$. Since 
$\support{H}$
is an intersecting family, $|J_u \intersect J_v| \geq 1$ and there exists some $j\in \Nset{m}$ with $j \in J_u \intersect J_v$. Thus, we have that $u, v\in c_j$ and since $c_j$ is a clique, $u \sim v$.  
\end{proof}
So a subgraph $H$ is connected if, and only if, its support is path-intersecting (Prop. \ref{prop:num_signatures_same_support}) and is a clique if, and only if, its support is an intersecting family (Prop. \ref{prop:cliques_to_intersecting_families}).
Theorem \ref{thm:maximal_clique} gives necessary and sufficient conditions for $H$ to be a \textit{maximal} clique.
%
\begin{theorem}
	\label{thm:maximal_clique}
	For any clique collection $\family{C} =\mCliqueSet{c}$,  a clique induced by $H$ on the graph union $\bigcup_{j =1}^{m} c_m$, is maximal, if, and only if,
	for any $J \subseteq \Nset{m}$,  
	\begin{enumerate}
		\item  
		 $J \in \support{H}$ 
		 $\imply$  $\cardinality{\Gamma_J \intersect H} = \gamma_J$, and
		\item 
		 $J \not \in \support{H}$ $\imply$ either $\Gamma_J = \emptyset$ or $\Gamma_J \ne \emptyset$ and $\{J\} \union \support{H}$ is not an intersecting family.
		 
		
	\end{enumerate}
\end{theorem}
\noindent
\begin{proof}

First, to prove necessity, assume $H$ is a maximal clique.  
For any $J \in \support{H}$, at least one node in $\Gamma_J$ is in $H$, and, so, connected to all other nodes in $H$.  It follows from Proposition \ref{prop:adjacency_intersection} that every node of $\Gamma_J$  is also in $H$ and hence $\cardinality{\Gamma_J \intersect H} = \gamma_J$ for all $J \in \support{H}$.  To show statement 2 holds, suppose now that $J \not \in \support{H}$.  Further, suppose that $\{J\} \union \support{H}$ is an intersecting family and so, by Proposition \ref{prop:cliques_to_intersecting_families}, that $\Gamma_J \union J$ is a clique.  Since $J \not \in \support{H}$, $\Gamma_J \intersect H = \emptyset$ and, since $H$ is maximal, it follows that $\Gamma_J = \emptyset$.

To prove sufficiency, assume $H$ is a clique and that both statements 1 and 2 hold.
By statement 1, all nodes in $\Gamma_J$ for $J \in \support{H}$ are in $H$ and no nodes remain in $\Gamma_J$ to increase $H$.  Statement 2 ensures that no nodes exist in any $\Gamma_J$ with $J \not \in \support{H}$ that could enlarge $H$ and still be a clique.
Hence, $H$ is maximal.
\end{proof}
Statement 2 of Theorem \ref{thm:maximal_clique} shows that, not only does a maximal clique have an intersecting family as its support (like all cliques), but that its intersecting family can only be expanded by sets $J \not \in \support{H}$ having no nodes in $\Gamma_J$.

\subsubsection{The $\Gamma$-quotient graph and maximal cliques}
\label{quotientandmaximalcliques}
Theorem  \ref{thm:maximal_clique} suggests that instead of considering intersecting families that are subsets of the entire power set, $\powerset{\Nset{m}}$, we need only those that are subsets of the support of the graph union $G = \union_{j=1}^{m}c_j$,  namely, $\support{G} = \{J : J\subseteq \Nset{m} \text{ and } \Gamma_J \neq \emptyset \} \subseteq \powerset{\Nset{m}}$.  

This effectively ignores empty cells of the $\Gamma$ partition to focus on intersecting families formed from the index sets that define the nodes of the quotient graph $G/\Gamma$.  The relevant families are intrinsic to the quotient graph.
For example, 
\begin{itemize}
\item any path on $G/\Gamma$ corresponds to a path-intersecting set (Prop. \ref{prop:num_signatures_same_support}),
\item any clique on $G/\Gamma$ determines an intersecting family and hence a clique on $G$, and 
\item any maximal clique on $G/\Gamma$ gives a maximal intersecting family and, so, a maximal clique on $G$.
\end{itemize}
The last two points are proved below in Proposition \ref{prop:maximal_cliques}.
\begin{proposition}
\label{prop:maximal_cliques}
	If $\mathcal{F}$ is a nonempty intersecting family on $\support{G}$, then the graph $H_\mathcal{F}$ induced by $\{\Gamma_J: J\in \mathcal{F}\}$ is a clique. Furthermore, $\mathcal{F}$ is a maximal intersecting family on $\support{G}$ if, and only if, $H_{\mathcal{F}}$ is a maximal clique.
\end{proposition}
\begin{proof}
The fact that is a clique follows immediately from Proposition \ref{prop:adjacency_intersection}. 

Suppose $\mathcal{F}$ is a maximal intersecting family on $\support{G}$ and $H_{\mathcal{F}}$ is not a maximal clique. Then there exists some $u \in V(G)$ with $u$ adjacent to all nodes in $H_\mathcal{F}$. Suppose $u\in \Gamma_{J_u}$, then $\Gamma_{J_u}$ is nonempty and by Proposition \ref{prop:adjacency_intersection}, $\Gamma_{J_u} \cap J \neq \emptyset$ for all $J \in \mathcal{F}$. Therefore, either $\mathcal{F}$ is not a maximal intersecting family or $H_{\mathcal{F}}$ was not the subgraph induced by $\mathcal{F}$ -- a contradiction.

The proof of the converse is almost identical.

\end{proof}
\begin{corollary}
\label{cor:unique_max_clique}
	If $\Gamma_J\neq \emptyset$ for all $\emptyset\neq J\subseteq \Nset{m}$, then every maximal intersecting family on $\powerset{\Nset{m}}$ induces a unique maximal clique in $G$.
\end{corollary}

\begin{proof}
Suppose $\Gamma_J\neq \emptyset$ for all $\emptyset\neq J\subseteq \Nset{m}$. Then $\support{G}$ is the set of all nonempty subsets of $\powerset{\Nset{m}}$. Therefore, by Proposition \ref{prop:maximal_cliques}, each maximal intersecting family gives to a unique maximal clique.
\end{proof}
This means that the number of maximal cliques, $M(\mathcal{C})$, in $G$ is equal to the number of maximal intersecting families on $\support{G}$ which in turn is bounded above by the number of maximal intersecting families on $\Nset{m}$.
\begin{corollary}
	\label{cor:maximal_clique_num}
The number, $M(\mathcal{C})$, of maximal cliques in the graph union of 
$\mathcal{C} = \mCliqueSet{c}$ is bounded above by $\lambda(m)$, the number of maximal intersecting families on $\Nset{m}$.
\end{corollary}
\begin{proof}
	By Theorem $\ref{thm:maximal_clique}$, each maximal intersecting family would correspond to at most one maximal clique in the graph union of the collection $\mCliqueSet{c}$. Thus, $\lambda(m)$ is an upperbound for $M(\mathcal{C})$. 
\end{proof}

\begin{corollary}
	\label{crl:clique_num}
	The clique number of the graph union of the collection of cliques $\mCliqueSet{c}$ is
\begin{align*}
\max_{\mathcal{F}\in \mathcal{M}} \quad & \sum_{J\in \mathcal{F}}\gamma_J,
\end{align*}
where $\mathcal{M}$ is the set of all maximal intersecting families on $\Nset{m}$.
\end{corollary}
\begin{proof}
By Theorem \ref{thm:maximal_clique}, a clique $H$ is maximal if, and only if, its corresponding intersecting family $\mathcal{F}_H$ is only extendible by trivial elements and $H$ uses all of the vertices in the cells $\Gamma_J$ that contain members from $H$. Therefore, for every maximal intersecting family $\mathcal{F}$, there is a corresponding unique maximal clique $H$ contained within the union of the cells $\{\Gamma_J: J\in F\}$. 

Since the clique number is the maximum of the size of all maximal cliques in a graph, and each maximal clique has the form $\sum_{J\in \mathcal{F}}\gamma_J$ for some maximal intersecting family $\mathcal{F}$, the proof follows.
\end{proof}
To summarize, an intersecting family on $\support{G}$ identifies a clique (Prop \ref{prop:cliques_to_intersecting_families}) and that clique is maximal if, and only if, its corresponding intersecting family is also maximal (Prop. \ref{prop:maximal_cliques}).  Whether an intersecting family, $\family{F}$, is maximally intersecting can be determined from its cardinality, namely an intersecting $\family{F} \subset \Nset{m}$ is a maximal intersecting family if, and only if, $|\mathcal{F}| = 2^{m-1}$ \cite<e.g., see Lemma 2.1>{meyerowitz1995maximal}; note that the intersecting family corresponding to an identified clique might have to be extended by adding subsets $J \in \Nset{m}$ having $\Gamma_J = \emptyset$ to achieve this cardinality (Thm. \ref{thm:maximal_clique}).  Every such maximal intersecting family produces a unique maximal
clique (Cor. \ref{cor:unique_max_clique}).  The number of such maximal cliques is bounded above by $\lambda(m)$, the number of maximally intersecting families on $\Nset{m}$ (Cor. \ref{cor:maximal_clique_num}). Unfortunately, $\lambda(m)$ is typically computationally intractable \cite<e.g., see>{brouwer2013counting} though is presently feasible on today's laptops for $m \le 10$, for example.
In the special case where $\gamma_J > 0$ for all $J\subseteq \Nset{m}$, every maximal intersecting family induces precisely one maximal clique so that the upper bound (Cor. \ref{cor:maximal_clique_num}) is achieved and $M(\mathcal{C}) = \lambda(m)$.

\section{Counting cliques}
\label{section:cliquecounting}
For a family of sets $\mathcal{F}$, let $N(\mathcal{F}) := \sum_{J\in \mathcal{F}}\gamma_J$ denote the number of nodes in the sets contained in the family. 

Given the collection of all maximal intersecting families on the support of $G$, we can apply the principle of inclusion-exclusion in the following manner.

\begin{proposition}
	\label{prop:cliques_containing_H}
	Let $H$ be a clique in the graph union of $\{c_1, \ldots, c_m\}$ and let $\mathcal{F}_H$ denote its support. Let $\mathcal{M}_H$ be the set of all maximal intersecting families $\mathcal{F}$ on $\support{G}$ that extend $\mathcal{F}_H$. The number of cliques that contain $H$ in the graph union of $\mCliqueSet{c}$ is
	\[ 
	1+\sum_{\mathcal{J}\subseteq \mathcal{M}_H}
			(-1)^{|\mathcal{J}|+1}
					\left(2^{N\left(\bigcap_{\mathcal{F} \in \mathcal{J}}
						\mathcal{F}\right) - |H|}-1\right).
	\]
\end{proposition}
\begin{proof}
Any clique that contains $H$ would be a subclique of one of the maximal cliques that contain $H$. Therefore, by Theorem $\ref{thm:maximal_clique}$, it suffices to examine the collection $\mathcal{M}_H$ of maximal intersecting families that generate a unique maximal clique in the graph union of $\mCliqueSet{c}$. If $\mathcal{F}\in \mathcal{M}_H$ corresponds to a maximal clique with $N(\mathcal{F})$ total nodes, then the selection of a nonempty subset from $\left(\cup_{J\in \mathcal{F}}\Gamma_J\right)\setminus H$ corresponds to a clique that properly contains $H$. This can be done in 
\[\left(2^{N\left(\mathcal{F}\right) - |H|}-1\right)\]
ways.

Since some cliques are are subgraphs of several different maximal cliques, we use the principle of inclusion/exclusion and obtain 
\[
\sum_{\mathcal{J}\subseteq \mathcal{M}_H}
			(-1)^{|\mathcal{J}|+1}
					\left(2^{N\left(\bigcap_{\mathcal{F} \in \mathcal{J}}
						\mathcal{F}\right) - |H|}-1\right)
	\]
cliques. However, this count does not include the clique $H$ on its own and hence we add a 1.
\end{proof}
The proof of Proposition \ref{prop:cliques_containing_H} relies on the fact that every clique is contained in some maximal clique.  
This observation can also be used to enumerate the total number of $r-$cliques in the graph union, by considering  
the collection of maximal cliques.

For instance, suppose we are interested in the number of triangles in the clique union from Figure $\ref{fig:three_K5_example}$.  There are only three maximal cliques in this graph union, each of size 5. 
To enumerate the triangles in the graph union, then, simply count the triangles in each maximal clique, subtract the number of triangles common to the each of the intersections of maximal cliques, and finally, add the number of triangles common to all three maximal cliques. This yields  
\begin{align*}
  \binom{5}{3} + \binom{5}{3} + \binom{5}{3} -\binom{2}{3} - \binom{3}{3} - \binom{3}{3} +\binom{2}{3} = 28
\end{align*}
triangles. 

The following proposition follows the same logic to generalize to counting the number of cliques of any size $r$ for any graph union of cliques.
Reminiscent of Proposition \ref{prop:PIE_for_cliques}, an advantage here is that the number of maximal cliques in the graph induced by the collection can be smaller than the number of cliques in the initial collection.
\begin{proposition}
	\label{prop:s_cliques_ext_cmb}
	The number of $r-$cliques induced by the graph union of the cliques $\mCliqueSet{c}$ is 
	\begin{align*}
  \sum_{\mathcal{J}\subseteq \mathcal{M}}(-1)^{|\mathcal{J}|+1}\binom{N(\intersect_{\mathcal{F}\in J}\mathcal{F})}{r},
\end{align*}
where $\mathcal{M}$ is the collection of all maximal intersecting families $\mathcal{F}$ with $\gamma_J > 0$ for all $J\in \mathcal{F}$.
\end{proposition}

\begin{proof}
As every clique is a subclique of a maximal clique, the induced graph by the maximal collection of cliques is the same as the induced graph by the collection $\mCliqueSet{c}$. Thus, the proof is exactly as in Proposition \ref{prop:PIE_for_cliques}.
\end{proof}

A more subtle expression for clique counts is had by considering their signatures.
\begin{theorem}
\label{thm:clique_gen_fn}
The generating function for clique counts induced by a collection $\mCliqueSet{c}$ is 
\[
\Phi(\sv{x}) =\sum_{\mathcal{F} \in \intersectingFamily{I}_m}
	\prod_{J \in \mathcal{F}}
		\left[(1+x_J)^{\gamma_J} - 1
		\right],
\]
where $\mathcal{I}$ is the set of all intersecting families on $\powerset{\Nset{m}}$, and $\sv{x}$ is the vector ${(x_J: J \in \powerset{\Nset{m}})}.$ 
\end{theorem}

\begin{proof}
A clique $H$ is determined uniquely by its signature and the node labels. By Proposition \ref{prop:cliques_to_intersecting_families}, the support must be an intersecting family on on $\support{G}$, and hence it is also an intersecting family on $\powerset{\Nset{m}}$. 

For a cell $J$ to contribute $\alpha_J\geq 1$ nodes to $H$ is accomplished in $\binom{\gamma_J}{\alpha_J}$ ways, which corresponds to the coefficient of $x_J^{\alpha_J}$ in the generating series 
\[ \left[(1+x_J)^{\gamma_J}-1 \right],\]
and the result follows.
\end{proof}
\noindent Extracting the coefficient of $x^r$ in the generating function $\Phi(x_J\to x)$ in Theorem \ref{thm:clique_gen_fn} yields the number of $r-$cliques as given in Corollary \ref{cor:r_cliques_by_signature}:
\begin{corollary}
	\label{cor:r_cliques_by_signature}
	The number of $r-$cliques in the graph union of the clique collection $\{c_1,\ldots c_m\}$ is 
	\[ 
	\sum_{\ell =1}^{r}
		\sum_{(\alpha_1,\ldots, \alpha_\ell)}
			\sum_{(J_1,\ldots, J_\ell)}
				\prod_{i=1}^\ell 
					\binom{\gamma_{J_i}}{\alpha_i}
	\] 
	where $(J_1,\ldots, J_\ell)$ is an intersecting family on $\support{G}$ of size $\ell$ with signature $(\alpha_1,\ldots,\alpha_\ell)$ being an integer composition of $r$ having $1 \leq \alpha_i\leq \gamma_i.$
\end{corollary}

\noindent 
A third expression for the total number of cliques of any size, induced by the collection, can also be had by substituting $x_J=1$ in the generating series in Theorem \ref{thm:clique_gen_fn}.  The expression is given as Corollary \ref{cor:total_clique_count}:
\begin{corollary}
\label{cor:total_clique_count}
The total number of cliques of size at least 1 induced by a collection $\mCliqueSet{c}$ is 
\[
\sum_{\mathcal{F} \in \intersectingFamily{I}_m}
	\prod_{J \in \mathcal{F}}
		\left[2^{\gamma_J} - 1
		\right],
\]
where $\intersectingFamily{I}_m$ is the set of all intersecting families on $\powerset{\Nset{m}}$.
\end{corollary}

\noindent When $r=2$, the interesting special case of the edge count is obtained (e.g., essential to edge count distributions for many random graph models, such as the Erd\H{o}s-R\'{e}nyi model):
\begin{corollary}
\label{cor:edge_count_by_type}
The number of edges induced by the collection of $r$-cliques $\mCliqueSet{c}$ is 

\begin{align*}
  \sum_{J \subseteq \Nset{m}} \binom{\gamma_J}{2} + \frac{1}{2} \sum_{J \subseteq \Nset{m}}\gamma_J 
  \sum_{I\neq J:~~|I\cap J| \geq 1}\gamma_I.
\end{align*}
\end{corollary}
Alternatively, edges can also be enumerated via
the degree sequences of the vertices in the various cells $\Gamma_J$. For every $J\subseteq \Nset{m}$, any two nodes within $\Gamma_J$ have the same degree. For instance, if $u \in \Gamma_{\{k\}}$ for some $k\in \Nset{m}$, then it must be that $deg(u) = r-1$ because $u \in c_k$ and $u\not\in c_j$ for all $j \neq k$ by the definition of $\Gamma_{\{k\}}$. On the other extreme, if $u\in \Gamma_{\Nset{m}}$, then $u \in c_j$ for all $j\in \Nset{m}$ and hence $u$ must be adjacent to all other nodes in $G$ which are in at least one of the $\{c_1, \ldots, c_m\}$. Therefore, $$deg(u) = n - \gamma_{\emptyset} - 1 = n - 1,$$
\noindent
The ``handshaking lemma'' immediately gives the number of edges induced by the collection as below:
\begin{proposition}
\label{prop:edge_count_by_deg}
The number of edges induced by the collection of cliques $\{c_1, \ldots, c_m\}$ is 
$$ \frac{1}{2}\sum_{J:\emptyset\neq J\subseteq \Nset{m}}\gamma_J
\left( \sum_{I: ~ |I\cap J| \geq 1}\gamma_I - 1\right).$$
\end{proposition}
\begin{proof}
	Follows immediately from Proposition \ref{prop:degree_seq_in_cliques} and the fact that number of edges in the graph is half the sum of the degrees in the graph.
\end{proof}


\section{Discussion}
\label{section:discussion}

In this work, connections were established and exploited between 
several graph-theoretic properties of clique covers, and notions of intersecting families on a special partition, the $\Gamma$-partition, of a graph $G = \union_{i=1}^m c_i$ formed from a collection of $\family{C} = \mCliqueSet{c}$ of $m$ cliques $c_i$. 

The partition was formed using elements $J$ of the power set $\powerset{\Nset{m}}$ from the
$m$ clique indices.  
The support of $G$ is a subset of the power set, $\support{G} \subseteq \powerset{\Nset{m}}$, and induces the partition $\Gamma$ of $\Nset{n}$, which partitions the set of $n$ distinct nodes in $G$.
This $\Gamma$-partition frames the unique contributions to $G$ from the various cliques of $\mCliqueSet{c}$ via sets from the power set of $\Nset{m}$. 


The quotient graph, $G/\Gamma$, induced by the $\Gamma$-partition succinctly captures the information provided by the collection of cliques.
This description serves as a dictionary between graph-theoretic traits, such as cliques, maximal cliques, and connected induced subgraphs, and their extremal set theory
counterparts  (viz., intersecting families, maximal intersecting families and path-intersecting families, respectively). 
The natural connection between these objects facilitates  
determination of expressions for several classes of counting problems arising from clique covers.

Of course, the $\Gamma$-partition and quotient graph are determined by the particular cliques given as elements of the collection.  Coarser partitions (those which produce fewer $\Gamma_J$ cells) are preferred -- ideally, the collection will consist of a minimal number of unique maximal cliques. 

Going forward, the techniques enabled by this partition approach may be adapted to enumerating graph components other than cliques (e.g., 
spanning trees or cycles).  
The methods also show promise in stochastic settings (e.g.,  since probability of particular graph configurations in Erd\H{o}s-R\'{e}nyi models is a function of edge counts, one can obtain the moments of clique counts on homogeneous 
Erd\H{o}s-R\'{e}nyi 
graphs using the techniques above). 

Finally, from a topological standpoint, we note that, since the number of $(r+1)-$cliques in a graph is corresponds to the number of $r$-faces in the clique complex of the graph, these results can be extended to study the bounds on the number of generators in the $r-$th homology class of the clique complex \cite<e.g.,>{kahle2009topology}. 
%
%

\bibliography{proposal_bibliography.bib}

\end{document}